\documentclass[11pt]{article}
             \newlength{\oldparindent}
             \usepackage{amsmath}
             \usepackage{latexsym}
             \usepackage{t1enc}
             \usepackage[latin1]{inputenc}

        \newtheorem{theoreme}{Theorem}[section]

             \newtheorem{lemme}[theoreme]{Lemma}
             \newtheorem{remarque}[theoreme]{Remark}
             \newtheorem{Not}[theoreme]{Notation}
             \newtheorem{hyp}[theoreme]{Assumption}

             \font\myf=msbm10 at 12pt
             \newcommand{\R}{\hbox{\myf R}}
             
             \newcommand{\pr}{\hbox{\myf P}}
             \newcommand{\esp}{\hbox{\myf E}}

             \newcommand{\F}{{\cal F}}

\begin{document}

  \def\second {\vtop{\baselineskip=11pt
             \hbox to 125truept{\hss 
             Diana DOROBANTU, Universit\'e de Lyon\footnote{Institut de Science Financi\`ere et d'Assurance, 50 avenue Tony Garnier, 69007 Lyon, France diana.dorobantu@adm.univ-lyon1.fr} \hss}\vskip2truept}}

            \title{{ A class of optimal stopping  problems for Markov processes  }}
             \author{ \second}
             \date{}
             \maketitle
        \textbf{Abstract} :  Our purpose is to study a particular class of optimal stopping  problems for Markov processes. 
We justify  the value function  convexity and  we deduce that there exists a boundary function  such that the   smallest  optimal stopping time is the first  time when the Markov process passes over the boundary depending on time. Moreover, we  propose a method to find the optimal boundary function.

 \vspace{0.3cm}
 
 \section{Introduction}
              In this paper we study  a particular optimal stopping problem for  strong Markov processes.  We propose a method to find the optimal stopping time form (it will be the first  time when the Markov process passes over a boundary depending on time), as well as for the calculation of the optimal boundary.

            In fact we seek to control a stochastic process V of the form $V=ve^X$ where $v$ is a real strictly positive constant and $X$ a  strong Markov process. We consider the following optimal stopping problem :
\begin{equation}
\nonumber
s(v)=sup_{ \tau\in \Delta} \esp\left[e^{-r\tau}h(V_{\tau}, ~\tau)\mid V_0=v\right],
\end{equation}
 where   $r>0$,  $\F^V_t=\sigma(V_s, ~s\leq t)$,  $\Delta$ is the set of   
$\F^V$-stopping times and $h$ is a Borelian function $h(V, t)=-V+ce^{mt}, ~c>0, ~m<r$.   We prove that our problem may be easily reduced to an optimal stopping problem for Markov processes and linear reward (i.e. $sup_{ \tau\in \Delta} \esp\left[e^{-r\tau}f(V_{\tau})\mid V_0=v\right]$ where $f$ is a linear function).   We justify  the convexity of the value function $s$ and  we deduce that the optimal strategy consists of stopping when the underlying Markov process crosses a boundary depending on time, i.e. the smallest optimal stopping time has the form $inf\{t\geq 0 : V_t\leq b(t)\}$. The main result is given by Theorems \ref{th1}, \ref{th2} and \ref{th3} which allow to determine the optimal stopping time form and the optimal boundary function. 

Optimal stopping theory  is a subject which often appears in the specialized literature.  For different areas of application or different methods for optimal stopping problems see, for example, Peskir and Shiryaev (2003). Among others, Salminen (1985), Leland (1994, 1996, 1998), Duffie and Lando (2001), Dayanik and Karatzas (2003) or Decamps and  Villeneuve (2007, 2008) studied  optimal stopping problems for  continuous Markov processes. Moreover, there are other authors who used  Lévy jumps processes  (e.g.   Pham (1997), Mordecki (1999), Hilberink and Rogers (2002), Kou and Wang (2004), Dao (2005), Kyprianou (2006), Dorobantu (2007)...) or symmetric Markov processes (e.g. Zabczyk (1984)) for their models. Sometimes the studied problem has the form $sup_{ \tau\geq 0} \esp\left[e^{-r\tau}h(V_{\tau})\right]$, other times it is more complicate $sup_{ \tau\geq 0} \esp\left[e^{-r\tau}h(V_{\tau}, \tau)\right].$ Our result completes these studies and the aim of the present paper is to solve a stopping time problem for a more general class of processes (more precisely, Markov processes not necessarily continuous). Contrary to the usual method, our method   avoids long calculations of the integro-differential operators.

This paper is organized as follows : we introduce the optimal stopping problem
(Section \ref{section2_ch1}). The following section (Section \ref{section3.1_ch1}) contains the main results which characterize the optimal stopping time and the optimal boundary.  Section \ref{appendices} is dedicated to  the proofs of Theorems \ref{th1}, \ref{th2} and \ref{th3}.
             \section{Optimal stopping problem}
              \label {section2_ch1}
              
                Let $V$ be a stochastic process on a filtered probability space $(\Omega,\F,(\F_t)_{t\geq 0},\pr)$. Assume that $V$ has the form  $V=ve^X$ where $v$ is  a real strictly positive constant and $X$ is a strong Markov process such that $X_0=0$.  Let $\F^V$ be  the right-continuous complete filtration generated by the process $V$, $\F^V_t=\sigma(V_s, ~s\leq t)$.  We introduce $\Delta$ the set of  $\F^V$-stopping times.

              From now on, $\esp(.|V_0=v)$ and $\pr(.|V_0=v)$ are denoted $\esp_v(.)$ and $\pr_v(.)$.

             We consider the following  optimal stopping  problem :
\begin {equation}
             \label{2}
            s(v)= sup_{\tau\in \Delta} \esp_v\left[e^{-r\tau}(-V_{\tau}+ce^{m\tau})\right],
             \end {equation}
             where  $r, ~c>0$ and $r>m$.

             We suppose that the process $X$ checks the following  assumptions  :
            \begin{hyp}
             \label{h0} 
             $\pr(lim_{t \downarrow 0} X_t=X_0)=1.$
             \end{hyp}
      \begin{hyp}
             \label{h1}
             The process $(e^{-rt+X_t}, ~t\geq 0)$ is of class $D$.
             \end{hyp}
             \begin{hyp}
             \label{h1i}
            $inf_{t\geq 0}e^{-rt}\esp(e^{X_t})=0$.
             \end{hyp}
             \begin{hyp}
             \label{h2}
            The support of $X_t$ is $\R$ for all $t>0$.
             \end{hyp}
             
            Under Assumptions \ref{h0}, \ref{h1}, \ref{h1i} and \ref{h2}, we  prove that the smallest optimal stopping time of (\ref{2}) is necessarily of the form $inf\{t\geq 0 : ~V_t\leq b(t)\}$ and we compute the optimal boundary function. We  applied the same method in \cite{d1, d2} for Lévy processes and linear functions (i.e. $m=0$), but it may be  extended to a more general class of processes and reward functions. The same type of problem as (\ref{2}) has been studied in \cite{d2} for a particular Markov process. The method used in \cite{d2} is different and it could be applied because the model is easy.

             \section{The main results}
            \label{section3.1_ch1}
            
            The main results caracterize the smallest optimal stopping time of (\ref{2}). We show the following.
           \begin{theoreme}
             \label{th1}
              Under Assumptions \ref{h0}, \ref{h1}, \ref{h1i} and \ref{h2}, there exists at least an optimal stopping time for the problem (\ref{2}).
                
                  For any  $c>0$, 
            there exists $b_c>0$  such that the smallest optimal stopping time has the following form 
             $$\tau_{b_c}=inf \{t\geq 0 : ~V_t \leq b_ce^{mt}\}.$$
             \end{theoreme}
             
              We introduce an auxiliary function 
           \begin {equation}
             \nonumber
             s_b(v)=\esp_v\left[e^{-(r-m)\tau_{b}}
             \left(- 
             e^{-m\tau_{b}}V_{\tau_{b}}+c\right)\right], ~~~v\in\R_+^*, ~~b\in \left]0, c\right[
        \end{equation}
        where $\tau_b=inf\{t\geq 0 : ~e^{-mt}V_t \leq b\}$. Let us point out that if $b\in\R_+$, then $s_b(.)$ is not necessarily positive. The condition $b\in ]0, c[$ implies the positivity of $s_b(.)$.
             \begin{remarque}
             \label{intervalle_ch1}
          Under  the assumptions of Theorem \ref{th1},  there exists $B_c$ such that $s_{B_c}(.)=s(.)$.
  \end{remarque}

Remark that we can write $s_.(.)$ as a function of Laplace transforms
\begin{equation}
              \nonumber
              {\cal{L}}(x)=\esp\left[e^{-(r-m)\bar{\tau}_x}|X_0=0\right],~~~  {\cal{G}}(x)=\esp\left[e^{-(r-m)\bar{\tau}_x+\bar{X}_{\bar{\tau}_x}}|X_0=0\right]
              \end{equation}
              where $\bar{X}$ is the process defined by $t\mapsto \bar{X}_t=-mt+X_t$ and $\bar{\tau}_x=inf\{t\geq 0 : \bar{X}_t\leq x\}$. Indeed, the function $s_.(.)$ can be written as
              $$s_b(v)=- v{\cal{G}}\left(ln\frac{b}{v}\right)+c{\cal{L}}\left(ln\frac{b}{v}\right).$$
              
              The following theorems caracterize the value of the optimal threshold $B_c$ as a function of  $c$, ${\cal{L}}(.)$ and ${\cal{G}}(.)$.
              
                When $\cal{G}$ is discontinuous at $x=0$, $B_c$ is easy to obtain.
               \begin{theoreme}
             \label{th2}
         Under Assumptions \ref{h0},  \ref{h1}, \ref{h1i} and \ref{h2}, we suppose that the function 
             $\cal{G}$  is discontinuous at $x=0$.  
             Then the smallest optimal stopping time is $\tau^*=inf\{t\geq 0 : V_t\leq 
             B_ce^{mt}\},$ where $B_c=c ~lim_{x\uparrow 0}\frac{1-{\cal{L}}(x)}{1-{\cal{G}}(x)}$.
             \end{theoreme}
             
              When $\cal{G}$ is continuous at $x=0$, $B_c$ is more technical to obtain, but it has the same form.
             \begin{theoreme}
             \label{th3}
             Under Assumptions \ref{h0},   \ref{h1}, \ref{h1i} and \ref{h2}, we suppose that the function 
             $\cal{G}$ is continuous at $x=0$. Then we have the following :
             \begin{enumerate}
\item
If $\cal{G}$ has  left derivative at $x=0$ (say  ${\cal{G}}'(0^-)$), then  $\cal{L}$ has left derivative at $x=0$ (say  ${\cal{L}}'(0^-)$). 
\item
 If moreover ${\cal{G}}'(0^-)\not=0,$ then $B_c\in [\tilde b , ~c[$ where $\tilde b=c  ~lim_{x\uparrow 0}\frac{1-{\cal{L}}(x)}{1-{\cal{G}}(x)}$.
 \item
If moreover $s_{\tilde b}(. )$ is strictly convex on $]\tilde b, ~\infty[$,
\end{enumerate}
     then the smallest optimal stopping time is $\tau^*=inf\{t\geq 0 : V_t\leq 
             B_ce^{mt}\}, \hbox{~where~} B_c=\tilde b.$
             \end{theoreme}
             
             The proofs of Theorems \ref{th1}, \ref{th2} and \ref{th3} are given in Section \ref{appendices}.
             
             \section{Appendix - Proofs}
             \label{appendices}
             
             Before starting with the proof of Theorem \ref{th1}, it is useful to re-formulate the problem (\ref{2}). For this purpose, following Gabillon (2003), we introduce a new  process $\nu$.
             \begin{Not}
             Let $\nu$ be the process defined by $\nu : t\mapsto ve^{-mt+X_t}(=ve^{\bar{X}_t})$.  We  sometimes use the notation $\nu^v=ve^{\bar{X}}$, for $v>0$. 
             \end{Not}

             The right-continuous complete filtration generated by the process $\nu$ is identical to $\F^V$. The problem (\ref{2}) may be written as 
             \begin{equation}
             \label{3}
            s(v)= sup_{ \tau\in \Delta} \esp_v\left[e^{-(r-m)\tau}f(\nu_{\tau})\right],
             \end{equation}
             where $f$ is a decreasing linear function, $f(v)=-v+c, ~v>0$. Therefore, problem (\ref{2}) can be reduced to an optimal stopping problem for Markov processes and linear functions.
             
             The proof of Theorem \ref{th1}  requires several results.

 Remark that  $s$ is a (decreasing) convex function because it is the sup of (decreasing) linear functions  :
             $$s(v)=sup_{\tau\geq 0}\esp_v\left[e^{-(r-m)\tau}(-
             \nu^v_{\tau}+c)\right]=sup_{\tau\geq 0}
             \esp_1\left[e^{-(r-m)\tau}(-  v\nu^1_{\tau}+c)\right].$$
             \begin{remarque}
           Since $s$ is a convex function, then it is   continuous.
             \end{remarque}

          The function $s$ is a positive function because
             $$s(v)\geq sup_{t\geq 0} \esp_v\left[e^{-(r-m)t} (- \nu_t+c)\right]\geq 
             sup_{t\geq 0} \esp_v\left[-e^{-(r-m)t}\nu_{t}\right]=sup_{t\geq 0}- v\esp\left[e^{-rt+X_t}\right]=0,$$
             where for the last equality we used Assumption \ref{h1i}.

 Under Assumption \ref{h1}, the process $\left(e^{-(r-m)t}f(\nu_t), ~t\geq 0\right)$ is of class D. According to Theorem 3.4 of \cite{KLM}, the Snell envelope of this process has the form $\left(e^{-(r-m)t}s(\nu_t), ~t\geq 0\right)$. Theorem  3.3 page 127   of \cite{Sh}, allows us to find the optimal stopping of a problem  $sup_{\tau\geq 0}\esp_v\left[f(\nu_{\tau})\right]$ where $f$ is a measurable function. We easily deduce  that this result may be applied to a process  having the form $t\mapsto e^{-rt}f(\nu_t)$. In our case, we can not apply this result for the problem (\ref{2}) because the process  $t\mapsto e^{-(r-m)t}f(\nu_t)$ does not check the assumptions of Theorem  3.3 page 127   of \cite{Sh} ; that is why we rewrite the function  $s$ under a new form. 
             \begin{lemme}
             \label{lem3}
               For $v>0$, let
             $s^+(v)=sup_{\tau\in\Delta}\esp_v\left[e^{-(r-m)\tau}(-\nu_{\tau}+c)^+\right],$
             where $x^+=max(x, 0)$.
Under Assumptions \ref{h0}, \ref{h1}, \ref{h1i} and \ref{h2},  $s^+(v)>0$ and $s(v)=s^+(v)$ for every $v>0$.
             \end{lemme}
             \begin{preuve}
             We show that if there exists $v_0>0$ such that  $s(v_0)<s^+(v_0)$, then there exists $v_1>0$ such that $s^+(v_1)=0$. We prove that this last relation  can not be satisfied.
             
           By construction, for each $v>0$, $s(v)\leq s^+(v)$.
             Let us suppose that there exists $v_0>0$ such that $s(v_0)<s^+(v_0)$.
             
      Under Assumption \ref{h0}, the process $\nu_.$ is  right continuous at 0. Since the process $Y^+ : t \rightarrow Y^+_t=e^{-(r-m)t}(-\nu_t+c)^+$ takes its values in  $[0, ~c]$, the assumptions of Theorem  3.3 page 127   of \cite{Sh} are checked for $Y^+$.
We denote by $f^+$ the function  $f^+(v)=(-\ v+c)^+$ ;  the stopping time
              $$\tau^+=inf\{u\geq 0 : f^+(\nu^{v_0}_u)=s^+(\nu^{v_0}_u)\}$$ is the smallest optimal stopping time of the problem $s^+(v_0)=sup_{\tau\geq 0}\esp_{v_0}\left[e^{-(r-m)\tau}(-\nu_{\tau}+c)^+\right].$

           Using the definition of  $s$ and $s^+$, we have
             $$\esp_{v_0}\left[ e^{-(r-m)\tau^+}f(\nu_{\tau^+})\right]\leq s(v_0)<s^+(v_0)=\esp_{v_0}\left[ e^{-(r-m)\tau^+}f^+(\nu_{\tau^+})\right]$$
             and consequently $$\esp_{v_0}\left[ e^{-(r-m)\tau^+}\left(f(\nu_{\tau^+})-f^+(\nu_{\tau^+})\right)\right]<0, ~~ 
             \pr_{v_0}\left(\{\omega : f(\nu_{\tau^+})<0\}\right)>0$$
            and $\pr_{v_0}\left(\{\omega : s^+(\nu_{\tau^+})=0\}\right)>0.$ 
             
            Thus there exists  $v_1$ such that $s^+(v_1)=0$.  Then for any stopping time $\tau$, $\pr_{v_1}$-almost surely $e^{-(r-m)\tau}f^+(\nu_{\tau})=0$ and  in particular for every  $t\in \R_+$, $f^+(\nu_t)=0$. This involves that $\pr_{v_1}$-almost surely $\nu_t\geq c$ which is a contradiction because under Assumption \ref{h2}, the support of $\nu_t$ is $\R^*_+$. Therefore $s^+(v)>0$ for every $v\in\R^*_+$ and $s(v)=s^+(v)$.
             \end {preuve}
             \hfill  $\Box$
             \\
             
             Thanks to  Lemma \ref{lem3}, the problem (\ref{2}) can be brought back to an optimal stopping problem for an American Put option with strike price $c$. Such a problem has been studied by many authors when $X$ is a Lévy process (see for exemple Gerber and Shiu (1994), Pham (1997), Mordecki (1999), Boyarchenko and Levendorskii (2002), Avram, Chan and Usabel (2002), Chesney and Jeanblanc (2004), Asmussen, Avram and Pistorius (2004), Alili and Kyprianou (2005), Kyprianou (2006)). Next,  we use a method close to the one used by Pham  (1997).  Pham  studies an optimal stopping problem for an American Put option with finite time horizon. In his model $X$ is a Lévy process. He uses integro-differential equations to solve his problem.
             \\
  \begin{preuve} \textbf{of Theorem \ref{th1}}
             \\
                       By Lemma \ref{lem3},  the problem  (\ref{2}) can be written as $sup_{\tau\geq 0}\esp(Y_{\tau }^+)$. By Theorem 3.3 page 127   of \cite{Sh},
               $\tau^*=inf\{u\geq 0 :
              f^+(\nu_u)=s^+(\nu_u)\}$ is the smallest optimal stopping time.
However $s(v)=s^+(v)>0$ for all $v>0$, so 
              $$\tau^*=inf\{u\geq 0 :
              f(\nu_u)=s(\nu_u)\}$$ is the smallest optimal stopping time.

                The function  $s$ is upper bounded by $c$  because $Y_.^+$ is upper bounded by
               $c$ and $lim_{v\downarrow 0} s(v)=lim_{v\downarrow 0}  f(v)=c$.
             
              Since  $s$ is convex, $f$ linear and $f(.)\leq s(.)$, then
             $\{v>0 : f(v)=s(v)\}$ is an interval of the form $]0, ~b_c]$.
    This means that the smallest optimal stopping time $\tau^*$ is also the first entrance time of $\nu$ in  
             $]0, ~b_c]$.
             \end {preuve}
             \hfill  $\Box$
             \\
             
             The smallest optimal stopping time is hence a hitting time for the process $\nu$.\\
            \begin{preuve} \textbf{of Theorem \ref{th2}}
             \\
             Let $b\in]0, c[$. The function $s_b(.)$ has the form 
             
                 \begin {equation}
             \nonumber
             s_b(v)=\left\{
             \begin {array}{l@{\quad}l}
             -v+ c& if \ v\leq b\\
        -v{\cal{G}}\left(ln\frac{b}{v}\right)+c{\cal{L}}\left(ln\frac{b}{v}\right) & if \ v>b. \\
             \end{array}
             \right.
             \end{equation}
             If the function $s_b(.)$ is continuous at $b$, then $b$ is solution of 
             \begin{equation}
             \label{bsoleq}
              -b+ c= -b{\cal{G}}(0^-)+c{\cal{L}}(0^-).
              \end{equation}
              However, $\cal{G}$ is discontinuous at $x=0$, so ${\cal{G}}(0^-)\not=1$ and the equation (\ref{bsoleq}) has only one solution : $$b^*=c\frac{1-{\cal{L}}(0^-)}{1-{\cal{G}}(0^-)}=c ~lim_{x\uparrow 0}\frac{1-{\cal{L}}(x)}{1-{\cal{G}}(x)}.$$
              
               The function $s$  has the form  $s_{B_c}(.)=s(.)$ and is convex, thus it is continuous, in particular it is continuous at $B_c$. We deduce that $B_c=b^*$.
           \end{preuve}
             \hfill  $\Box$
             \\
            \begin{preuve} \textbf{of Theorem \ref{th3}}
             \\
             $(1)$ By Remark \ref{intervalle_ch1},  there exists $B_c$ such that $s_{B_c}(.)=s(.)$. The function $s$ is convex,  therefore the right and  left derivatives exist everywhere and
             \begin{equation}
             \label{deriv_s}
             s'(v^-)\leq s'(v^+) \hbox{~for all~} v\in\R^*_+,
             \end{equation}
             where   $s'(v^-)$ and $s'(v^+)$ are the left and right derivatives of $s$ at $v$. In particular, this means that
             $$s_{B_c}(v)=-v{\cal{G}}\left(ln\frac{B_c}{v}\right)+c{\cal{L}}\left(ln\frac{B_c}{v}\right)=s(v)$$
             has right and  left derivatives at $v=B_c$.
             Since $\cal{G}$ has right and left derivatives at $x=0$, then $\cal{L}$ has also right and left derivatives at $x=0$.
             \\
             $(2)$ Let us make $v=B_c$ in (\ref{deriv_s}) :
             $$-1\leq -1+{\cal{G}}'(0^-)-\frac{c}{B_c}{\cal{L}}'(0^-).$$
                         We deduce that  $B_c\geq\tilde b=c ~\frac{{\cal{L}}'(0^-)}{{\cal{G}}'(0^-)}=c ~lim_{x\uparrow 0}\frac{1-{\cal{L}}(x)}{1-{\cal{G}}(x)}$.
   \\
       $(3)$  If moreover    $s_{\tilde b}(.)$ is strictly convex on $]\tilde b, ~\infty[$, then
       \begin{equation}
       \label{gf}
s_{\tilde b}(v)>f(v) \hbox{~for all~} v> \tilde b.
\end{equation}
       Indeed, the graph of  $f$ is tangent to the graph of $s_{\tilde b}(.)$ in $v=\tilde b$. 
       
       Suppose that $B_c> \tilde b$, then 
       $f(B_c)=s(B_c)=s_{B_c}(B_c)\geq s_{\tilde b}(B_c)$   which contradicts    (\ref{gf}).    
              \end{preuve}
             \hfill  $\Box$
             \\
             \begin{remarque}
             Assumption \ref{h1} may be replaced by 
  
           "There exists $q\in\R$ such that the support of $X_t$ is included in $]-\infty, ~q]$ for all $t>0$."

           Under this assumption, we don't need to use the intermediate Lemma \ref{lem3} to find the smallest optimal stopping time form. In this case the process $(f(\nu_t), ~t\geq 0)$ is bounded and Theorem 3.3 page 127   of \cite{Sh} can be  directly applied. The function $s$ is not necessarily continuous, but its continuous extension by linear interpolation is convex and the conclusion of Theorems \ref{th1}, \ref{th2} and \ref{th3} are true.
           \end{remarque}
           
          Our results are consistent with existing literature. Recall that our problem can be brought back to an American Put optimal stopping problem for strong Markov processes. Various authors have found that, in the case of a Lévy process, the American Put optimal stopping problem is linked to the first passage problem of the Lévy process. Moreover, the optimal threshold is obtained using continuous or smooth pasting condition. For example, in \cite{Ali, boy} sufficient or necessary and sufficient conditions for smooth and continuous pasting were established for different classes of Lévy processes. To this subject (but for a different optimal stopping problem), see also \cite{KypS}. The aim of this paper is to solve a little more general  problem than the American Put optimal stopping problem, for a more general class of processes.

             \end{document}